\documentclass[10pt, a4paper,  reqno]{amsart}
\usepackage{amsfonts}
\usepackage{amssymb}
\usepackage{amsmath}

\newtheorem{theorem}{Theorem}
\newtheorem{lemma}[theorem]{Lemma}
\newtheorem{corollary}[theorem]{Corollary}

\newcommand{\cx}{\mathbb{C}}
\newcommand{\nt}{\mathbb{N}}
\newcommand{\iz}{\mathbb{Z}}
\newcommand{\ra}{\mathbb{Q}}
\newcommand{\fps}{\mathcal{F}}
\newcommand{\pp}{\mathbb{P}}
\newcommand{\B}[2]{B_{#1}(#2)}
\newcommand{\Br}[3]{B_{#1}^{(#2)}(#3)}

\newcommand{\vs}{\vspace{2mm}}

\begin{document}

\title[Umbral Calculus associated with Bernoulli polynomials]{Umbral Calculus associated with Bernoulli polynomials}

\author{Dae San Kim$^{1}$}
\address{$^{1}$ Department of Mathematics \\
Sogang University, Seoul 121-741, Republic of Korea  \\}
\email{dskim@sogang.ac.kr}

\author{Taekyun Kim$^{2}$}
\address{$^{2}$ Department of Mathematics \\
Kwangwoon University, Seoul 139-701, Republic of Korea  \\}
\email{tkkim@kw.ac.kr}

\begin{abstract}
Recently, D. S. Kim  and T. Kim have studied applications of umbral calculus associated with $p$-adic invariant integrals on $\Bbb Z_p$
(see [6]). In this paper, we
investigate some interesting properties arising from umbral
calculus. These properties are useful in deriving some identities of
Bernoulli polynomials.
\end{abstract}


\maketitle

\section{Introduction}

As is well known, the Bernoulli polynomials are defined by the
generating function to be
\begin{equation}\label{eqn_01}
\frac{te^{xt}}{{e^t}-1} = e^{B(x)t} = \sum_{n=0}^{\infty}
B_n(x)\frac{t^n}{n!},
\end{equation}
with the usual convention about replacing $B^n(x)$ by $B_n(x)$. In
the special case, $x=0$, $B_n(0)= B_n$ are called the $n$-th
Bernoulli numbers. From $(\ref{eqn_01})$, we note that
\begin{equation*}
B_0 = 1, \quad (B+1)^n- B^n = B_n(1) - B_n = \delta_{1,n},\quad
(\text{see \cite{02,03,04}}),
\end{equation*}
where $\delta_{m,k}$ is the Kronecker symbol.

\vs

In particular, by (\ref{eqn_01}), we set
\begin{equation}\label{eqn_02}
B_n(x) = \sum_{l=0}^n \binom{n}{l}B_{n-l}x^l,\quad (\text{see
\cite{01,02,03,04,06,07}}).
\end{equation}
By (\ref{eqn_02}), we see that $B_n(x)$ is a monic polynomial of
degree $n$. We recall the Euler polynomials are defined by the
generating function to be
\begin{equation}\label{eqn_03}
\frac{2e^{xt}}{{e^t}+1} = e^{E(x)t} = \sum_{n=0}^{\infty}
E_n(x)\frac{t^n}{n!} ,\quad (\text{see \cite{01,07}}),
\end{equation}
with the usual convention about replacing $E^n(x)$ by $E_n(x)$. In
the special case, $x=0$, $E_n(0)= E_n$ are called the $n$-th Euler
numbers. From $(\ref{eqn_03})$, we can derive the following
equation:
\begin{equation}\label{eqn_04}
E_n(x) = \sum_{l=0}^n \binom{n}{l}E_{n-l}\ x^l,\quad (\text{see
\cite{07,09}}).
\end{equation}
Thus (\ref{eqn_04}), we see that $E_n(x)$ is also a monic polynomial of degree $n$.
By (\ref{eqn_04}), we get
\begin{equation}\label{eqn_05}
E_0 = 1, \quad (E+1)^n+ E^n = E_n(1) + E_n = 2\ \delta_{0,n}.
\end{equation}

\vs

Let $\cx$ be the complex number field and let $\fps$ be the set of
all formal power series in the variable $t$ over $\cx$ with
\begin{equation*}
\fps = \left\{ f(t) = \sum_{k=0}^\infty \frac{a_k}{k!} t^k : a_k \in
\cx \right\}.
\end{equation*}
We use the notation $\pp = \cx[x]$ and $\pp^*$
denotes the vector space of all linear functional on $\pp$.

\vs

Let $\left\langle L\ | \ p(x)\right\rangle $ be the action of a linear functional $L$ on a
polynomial $p(x)$, and we remaind that the vector space operation on
$\pp^*$ are defined by
\begin{equation*}
\begin{split}
\left\langle L+M\ | \ p(x)\right\rangle  &= \left\langle L\ | \ p(x)\right\rangle  + \left\langle M\ | \ p(x)\right\rangle , \\
\left\langle cL\ | \ p(x)\right\rangle  &= c\left\langle L\ | \ p(x)\right\rangle , \quad (\text{see \cite{05,16}}),
\end{split}
\end{equation*}
where $c$ is any constant in $\cx$.

\vs

The formal power series
\begin{equation}\label{eqn_06}
f(t) = \sum_{k=0}^{\infty} \frac{a_k}{k!} t^k \in \fps
\end{equation}
defines a linear functional on $\pp$ by setting
\begin{equation}\label{eqn_07}
\left\langle f(t)\ | \ x^n\right\rangle  = a_n, \ \text{for all }n\geq 0.
\end{equation}
Thus, by (\ref{eqn_06}) and (\ref{eqn_07}), we have
\begin{equation}\label{eqn_08}
\left\langle t^k\ | \ x^n\right\rangle  = n! \ \delta_{n,k}, \quad (\text{see \cite{05,16}}).
\end{equation}

\vs

Let $f_L(t) = \sum_{k=0}^{\infty} \frac{\left\langle L\ | \ x^k\right\rangle }{k!} t^k$. Then we
see that $\left\langle f_L(t)\ | \ x^n\right\rangle  = \left\langle L\ | \ x^n\right\rangle $ and so as linear functionals
$L=f_L(t)$ (see \cite{05,16}). As is known in \cite{16}, the map $L
\mapsto f_L(t)$ is a vector space isomorphism from $\pp^*$ onto
$\fps$. Henceforth, $\fps$ will denote both the algebra of formal
poser series in $t$ and the vector space of all linear functionals
on $\pp$, and so an element $f(t)$ of $\fps$ will be thought of as
both a formal power series and a linear functional. We shall call
$\fps$ the umbral algebra. The umbral calculus is the study of
umbral algebra and modern classical umbral calculus can be described
as a systematic study of the class of Sheffer sequences (see
\cite{16}).

\vs

The order $ord(f(t))$ of a nonzero power series $f(t)$ is the smallest
integer $k$ for which the coefficient of $t^k$ does not vanish. If a series $f(t)$ with $ord(f(t))=1$, then $f(t)$ is called a delta
series. If a series $f(t)$ with $ord(f(t))=0$, then $f(t)$ is called
an invertible series (see \cite{05,16}). For $f(t),g(t)\in \fps$, we
have
\begin{equation}\label{eqn_09}
\left\langle f(t)g(t)\ | \ p(x)\right\rangle  = \left\langle f(t) \ | \  g(t)p(x)\right\rangle , \ (\text{see \cite{05,16}}).
\end{equation}

\vs

Let us assume that $S_n(x)$ denotes a polynomial of degree $n$. If
$f(t)$ is a delta series and $g(t)$ is an invertible series, then there
exists a unique sequence $S_n(x)$ such that $\left\langle g(t)f(t)^k \ | \  S_n(x)\right\rangle  =
n! \ \delta_{n.k}$, $n,k\geq 0$ (see \cite{05}). The sequence $S_n(x)$
is called the Sheffer sequence for $(g(t), f(t))$, denoted by $S_n(x)
\sim(g(t), f(t))$. If $S_n(x) \sim(1, f(t))$, then $S_n(x)$ is called the associated
sequence for $f(t)$ or $S_n(x)$ is associated to $f(t)$. If $S_n(x)
\sim(g(t), t)$, then $S_n(x)$ is called the Appell sequence for
$g(t)$ or $S_n(x)$ is Appell for $g(t)$ (see \cite{05,16}). For
$p(x)\in \pp$, it is known (see \cite{05,16}) that
\begin{equation}\label{eqn_10}
\left\langle \frac{e^{yt}-1}{t}\ | \ p(x)\right\rangle  = \int_{0}^{y} p(u) du,
\end{equation}
\begin{equation}\label{eqn_11}
\left\langle f(t)\ | \ xp(x)\right\rangle  = \left\langle \partial_tf(t)\ | \ p(x)\right\rangle  = \left\langle f^{'}(t)\ | \ p(x)\right\rangle ,
\end{equation}
and
\begin{equation}\label{eqn_12}
\left\langle e^{yt} - 1\ | \ p(x)\right\rangle  = p(y) - p(0), \ (\text{see \cite{05,16}}).
\end{equation}

\vs

Let us assume that $S_n(x) \sim(g(t), f(t))$. Then we have the
following equations (\ref{eqn_13})-(\ref{eqn_16}):
\begin{equation}\label{eqn_13}
h(t) = \sum_{k=0}^{\infty} \frac{\left\langle h(t)\ | \ S_k(x)\right\rangle }{k!} g(t) f(t)^k, \
h(t) \in\fps,
\end{equation}
\begin{equation}\label{eqn_14}
p(t) = \sum_{k=0}^{\infty} \frac{\left\langle g(t)f(t)^k\ | \ p(x)\right\rangle }{k!} S_k(x), \
p(t) \in\pp,
\end{equation}
\begin{equation}\label{eqn_15}
f(t)S_n(x) = n S_{n-1}(x), \ (n \in \iz_{+} = \nt \cup \{0 \}),
\end{equation}
and
\begin{equation}\label{eqn_16}
\frac{1}{g(\bar{f}(t))}e^{y \bar{f}(t)} = \sum_{k=0}^{\infty}
\frac{S_k(y)}{k!} t^k, \ \text{ for all }y \in \cx.
\end{equation}
where $\bar{f}(t)$ is the compositional inverse of $f(t)$ (see
\cite{16}).

\vs

Let $f_1(t), \ldots, f_m(t) \in \fps$. Then as is well known, we
have
\begin{equation}\label{eqn_17}
\left\langle f_1(t)f_2(t)\cdots f_m(t)\ | \ x^n\right\rangle  = \sum \binom{n}{i_1, \ldots, i_m}
\left\langle f_1(t) \ | \  x^{i_1}\right\rangle  \cdots \left\langle f_m(t) \ | \  x^{i_m}\right\rangle ,
\end{equation}
where the sum is over all nonnegative integers $i_1, \ldots, i_m$ such that $i_1 + \cdots
+ i_m = n$ (see \cite{05,16}).

\vs

In \cite{05}, D. S. Kim and T. Kim have studied  applications of umbral calculus associated with $p$-adic invariant integrals on $\Bbb Z_p$.
In this paper, we derive some
interesting properties of Bernoulli polynomials arising from umbral
calculus. These properties will be used in studying identities on the
Bernoulli polynomials

\section{Umbral Calculus and Bernoulli polynomials}

Let $\pp_n=\{p(x) \in \cx[x] : \deg p(x) \leq n\}$ and let $S_n(x)
\sim (g(t),t)$. From (\ref{eqn_16}), we have
\begin{equation}\label{eqn_18}
\frac{1}{g(t)} x^n = S_n(x) \ \Leftrightarrow \ x^n = g(t)S_n(x), \
(n\geq 0).
\end{equation}
Let us take $g(t) = \frac{1}{t}(e^t - 1) \in \fps$. Then $g(t)$ is
invertible series. By (\ref{eqn_01}), we get
\begin{equation}\label{eqn_19}
\sum_{k=0}^{\infty} \frac{B_k(x)}{k!} t^k = \frac{1}{g(t)} e^{xt}.
\end{equation}
Thus by (\ref{eqn_19}), we have
\begin{equation}\label{eqn_20}
\frac{1}{g(t)} x^n = B_n(x), \ (n\geq 0),
\end{equation}
and
\begin{equation}\label{eqn_21}
tB_n(x) = B_n'(x) = n B_{n-1}(x).
\end{equation}
From (\ref{eqn_20}) and (\ref{eqn_21}), we note that $B_n(x)$ is an
Appell sequence for $\frac{1}{t}(e^t - 1)$. By (\ref{eqn_02}), we
get
\begin{equation}\label{eqn_22}
\begin{split}
\int_{x}^{x+y} B_n(u)du & = \frac{1}{n+1} \{B_{n+1}(x+y) -
B_{n+1}(x)\} \\
& = \sum_{k=1}^{\infty}\frac{y^k}{k!} t^{k-1}B_n(x) =
\frac{e^{yt}-1}{t} B_n(x).
\end{split}
\end{equation}
In particular, for $y=1$, we have
\begin{equation}\label{eqn_23}
B_n(x) = \frac{t}{e^{t}-1} \int_{x}^{x+1} B_n(u)du = \frac{t}{e^t -
1} x^n.
\end{equation}
By (\ref{eqn_15}), we easily get
\begin{equation}\label{eqn_24}
B_n(x) = t \left\{\frac{1}{n+1} B_{n+1}(x) \right\}.
\end{equation}
From (\ref{eqn_24}), we can derive the following equation:
\begin{equation}\label{eqn_25}
\left\langle \frac{e^{yt}-1}{t}\ | \ B_n(x)\right\rangle  = \left\langle e^{yt}-1 \ | \  \frac{1}{n+1}B_{n+1}(x)\right\rangle  =
\int_{0}^{y} B_n(u) du.
\end{equation}

\vs

For $r\in \nt$, the $n$-th Bernoulli polynomials of order $r$ are
defined by the generating function to be
\begin{equation}\label{eqn_26}
\begin{split}
\left(\frac{t}{e^{t}-1}\right)^r e^{xt}
& =
\sum_{n=0}^{\infty} B_n^{(r)}(x) \frac{t^n}{n!}, \ (\text{see
\cite{02,03}}).
\end{split}
\end{equation}

\vs

In the special case, $x=0$, $B_n^{(r)}(0) = B_n^{(r)}$ are called
the $n$-th Bernoulli numbers of order $r$. By (\ref{eqn_26}), we get
\begin{equation}\label{eqn_27}
B_n^{(r)}(x) = \sum_{l=0}^{n} \binom{n}{l} B_{n-l}^{(r)} \ x^r, \
(\text{see \cite{02,03}}).
\end{equation}
Note that
\begin{equation}\label{eqn_28}
B_n^{(r)} = \sum_{l_1 + \cdots + l_r = n} \binom{n}{l_1,\ldots, l_r}
B_{l_1}\cdots B_{l_r}.
\end{equation}
From (\ref{eqn_27}) and (\ref{eqn_28}), we note that $B_n^{(r)}(x)$
is a monic polynomials with coefficients in $\ra$. By
(\ref{eqn_27}), we get
\begin{equation}\label{eqn_29}
\begin{split}
\int_{x}^{x+y} \Br{n}{r}{u} du & = \frac{1}{n+1} \{\Br{n+1}{r}{x+y} - \Br{n+1}{r}{x}\} \\
& = \sum_{k=1}^{\infty} \frac{y^k}{k!} t^{k-1} \B{n}{x} = \frac{e^{yt}-1}{t} \Br{n}{r}{x},
\end{split}
\end{equation}
and
\begin{equation}\label{eqn_30}
\Br{n}{r}{x+1} - \Br{n}{r}{x} = n \Br{n-1}{r-1}{x}.
\end{equation}
From (\ref{eqn_29}) and (\ref{eqn_30}), we note that
\begin{equation}\label{eqn_31}
\frac{e^t-1}{t}\Br{n}{r}{x} = \int_{x}^{x+1} \Br{n}{r}{u} du = \Br{n}{r-1}{x}.
\end{equation}
By (\ref{eqn_31}), we get
\begin{equation}\label{eqn_32}
\Br{n}{r}{x} = \left(\frac{t}{e^t-1}\right) \Br{n}{r-1}{x} = \left(\frac{t}{e^t-1}\right)^{r-1} \B{n}{x} = \left(\frac{t}{e^t-1}\right)^{r} x^n,
\end{equation}
and
\begin{equation}\label{eqn_33}
t\Br{n}{r}{x} = n \left(\frac{t}{e^t-1}\right)^{r} x^{n-1} = n \Br{n-1}{r}{x}.
\end{equation}
It is easy to show that $\left( \frac{e^t-1}{t}\right)^r$ is an invertible series in $\fps$. Therefore, by (\ref{eqn_32}) and (\ref{eqn_33}), we obtain the following lemma.
\begin{lemma}
$\Br{n}{r}{x}$ is the Appell sequence for $\left( \frac{e^t-1}{t}\right)^r$.
\end{lemma}
By (\ref{eqn_33}), we get
\begin{equation}\label{eqn_34}
\Br{n}{r}{x} = t \left\{ \frac{1}{n+1} \Br{n+1}{r}{x}\right\}, \ (n\geq 0).
\end{equation}
Thus, from (\ref{eqn_34}), we have
\begin{equation}\label{eqn_35}
\left\langle \frac{e^{yt}-1}{t}\ | \ \Br{n}{r}{x}\right\rangle  = \left\langle e^{yt}-1\ | \ \frac{1}{n+1}\Br{n+1}{r}{x}\right\rangle  = \int_{0}^{y} \Br{n}{r}{u} du.
\end{equation}

\vs

In the special case, $y=1$, we have
\begin{equation*}
\left\langle \frac{e^{t}-1}{t}\ | \ \Br{n}{r}{x}\right\rangle  = \int_{0}^{1} \Br{n}{r}{u} du = B_{n}^{(r-1)}.
\end{equation*}
By (\ref{eqn_17}), we get
\begin{equation}\label{eqn_36}
\small
\left\langle \left(\frac{t}{e^{t}-1}\right)^r\ | \ x^n\right\rangle  = \sum_{n=i_1+\cdots+i_r} \binom{n}{i_1,\ldots, i_r} \left\langle \frac{t}{e^{t}-1}\ | \ x^{i_1}\right\rangle  \cdots \left\langle \frac{t}{e^{t}-1}\ | \ x^{i_r}\right\rangle
\end{equation}
and
\begin{equation}\label{eqn_37}
\left\langle \frac{t}{e^{t}-1}\ | \ x^n\right\rangle  = B_n, \quad \left\langle \left(\frac{t}{e^{t}-1}\right)^r\ | \ x^{n}\right\rangle  = B_n^{(r)}.
\end{equation}
Thus, from (\ref{eqn_36}) and (\ref{eqn_37}), we have
\begin{equation*}
\sum_{n=i_1+\cdot + i_r} \binom{n}{i_!,\ldots,i_r} B_{i_1}\cdots B_{i_r} = B_n^{(r)}.
\end{equation*}

\vs

Let us take $p(x) \in \pp_n$ with
\begin{equation}\label{eqn_38}
p(x)= \sum_{k=0}^{n} b_k \B{k}{x}.
\end{equation}
From (\ref{eqn_20}) and (\ref{eqn_21}), we note that $\B{n}{x} \sim (\frac{e^{t}-1}{t},t)$. By the definition of Appell sequences, we get
\begin{equation}\label{eqn_39}
\left\langle \frac{e^{t}-1}{t}t^k\ | \ \B{n}{x}\right\rangle  = n! \ \delta_{n,k}, \ (n,k\geq 0),
\end{equation}
and, from (\ref{eqn_38}), we have
\begin{equation}\label{eqn_40}
\left\langle \frac{e^{t}-1}{t}t^k\ | \ p(x)\right\rangle  = \sum_{l=0}^{n} b_l \left\langle \frac{e^{t}-1}{t}t^k\ | \ \B{l}{x}\right\rangle  = \sum_{l=0}^{n} b_l l! \delta_{l,k} = k! b_k.
\end{equation}
Thus, by (\ref{eqn_25}) and (\ref{eqn_40}), we get
\begin{equation}\label{eqn_41}
b_k = \frac{1}{k!} \left\langle \frac{e^{t}-1}{t}t^k\ | \ p(x)\right\rangle  = \frac{1}{k!} \left\langle \frac{e^{t}-1}{t}\ | \ p^{(k)}(x)\right\rangle  = \frac{1}{k!} \int_{0}^{1} p^{(k)}(u) du,
\end{equation}
where $p^{(k)}(u) = \frac{d^k}{du^k}p(u)$. Therefore, by (\ref{eqn_38}) and (\ref{eqn_41}), we obtain the following theorem.

\begin{theorem}\label{thm_02}
Let $p(x)\in \pp_n$ with $p(x) = \sum_{k=0}^{n} b_k \B{k}{x}$. Then we have
\begin{equation*}
b_k = \frac{1}{k!} \left\langle \frac{e^{t}-1}{t}\ | \ p^{(k)}(x)\right\rangle  = \frac{1}{k!} \int_{0}^{1} p^{(k)}(u) du,
\end{equation*}
where $p^{(k)}(u) = \frac{d^k}{du^k}p(u)$.
\end{theorem}

\vs

Let $p(x) = \Br{n}{r}{x} \in \pp_n$ with $p(x) = \sum_{k=0}^{n} b_k \B{k}{x}$. Then we have
\begin{equation}\label{eqn_42}
p^{(k)}(x) = k! \binom{n}{k} \Br{n-k}{r}{x},
\end{equation}
and
\begin{equation}\label{eqn_43}
b_k = \frac{1}{k!} \left\langle \frac{e^{t}-1}{t}\ | \ p^{(k)}(x)\right\rangle  = \binom{n}{k} \left\langle \frac{e^{t}-1}{t}\ | \ \Br{n-k}{r}{x}\right\rangle .
\end{equation}
Therefore, by Theorem \ref{thm_02} and (\ref{eqn_43}), we obtain the following theorem.
\begin{corollary}
For $n\geq 0$, we have
\begin{equation*}
\Br{n}{r}{x} = \sum_{k=0}^{n}\binom{n}{k} \left\langle \frac{e^{t}-1}{t}\ | \ \Br{n-k}{r}{x}\right\rangle  \B{k}{x}.
\end{equation*}
In other words,
\begin{equation*}
\Br{n}{r}{x} = \sum_{k=0}^{n}\binom{n}{k} B_{n-k}^{(r-1)}\B{k}{x}.
\end{equation*}
\end{corollary}

\vs

From the definition of Appell sequences, we note that
\begin{equation}\label{eqn_45}
\left\langle \left(\frac{e^{t}-1}{t}\right)^r t^k \ | \  \Br{n}{r}{x}\right\rangle  = n! \ \delta_{n,k}, \ (n,k\geq 0).
\end{equation}

\vs

Let $p(x) \in \pp_n$ with $p(x) = \sum_{k=0}^{n} b_k^{(r)} \Br{k}{r}{x}$. By (\ref{eqn_45}), we get
\begin{equation}\label{eqn_47}
\begin{split}
\left\langle \left(\frac{e^{t}-1}{t}\right)^r t^k \ | \  p(x) \right\rangle  & = \sum_{l=0}^{n} b_l^{(r)} \left\langle \left(\frac{e^{t}-1}{t}\right)^r t^k \ | \  \Br{l}{r}{x} \right\rangle  \\
& = \sum_{l=0}^{n} b_l^{(r)} l! \ \delta_{l,k} = k! \ b_k^{(r)}.
\end{split}
\end{equation}
Thus, by (\ref{eqn_47}), we have
\begin{equation}\label{eqn_48}
b_k^{(r)} = \frac{1}{k!} \left\langle \left(\frac{e^{t}-1}{t}\right)^r t^k \ | \  p(x)\right\rangle .
\end{equation}
Therefore, by (\ref{eqn_48}), we obtain the following theorem.
\begin{theorem}\label{thm_03}
Let $p(x) \in \pp_n$ with $p(x) = \sum_{k=0}^{n} b_k^{(r)} \Br{k}{r}{x}$. Then we have
\begin{equation*}
b_k^{(r)} = \frac{1}{k!} \left\langle \left(\frac{e^{t}-1}{t}\right)^r t^k \ | \  p(x)\right\rangle .
\end{equation*}
\end{theorem}

\vs

Let us consider $p(x) = \B{n}{x}$ with
\begin{equation}\label{eqn_49}
\B{n}{x} = p(x) = \sum_{k=0}^{n} b_k^{(r)} \Br{k}{r}{x}.
\end{equation}
By Theorem \ref{thm_03} and (\ref{eqn_49}), we get
\begin{equation}\label{eqn_50}
b_k^{(r)} = \frac{1}{k!} \left\langle \left(\frac{e^{t}-1}{t}\right)^r t^k \ | \  p(x)\right\rangle    = \frac{1}{k!} \left\langle \left(\frac{e^{t}-1}{t}\right)^r t^k \ | \  \B{n}{x}\right\rangle .
\end{equation}

\vs

For $k< r$, we have
\begin{equation}\label{eqn_51}
\begin{split}
b_k^{(r)} & = \frac{1}{k!} \left\langle \left(\frac{e^{t}-1}{t}\right)^r t^k \ | \  \B{n}{x}\right\rangle  \\
 & = \frac{1}{k!} \left\langle (e^{t}-1)^r \ | \  \frac{\B{n+r-k}{x}}{(n+1) \cdots (n+r-k)}\right\rangle  \\
 & = \frac{1}{k!(r-k)!\binom{n+r-k}{r-k}} \left\langle (e^{t}-1)^r \ |  \B{n+r-k}{x}\right\rangle  \\
 & = \frac{\binom{r}{k}}{r!\binom{n+r-k}{r-k}} \sum_{j=0}^{r} \binom{r}{j} (-1)^{r-j} \left\langle e^{jt} \ | \  \B{n+r-k}{x}\right\rangle  \\
 & = \frac{\binom{r}{k}}{r!\binom{n+r-k}{r-k}} \sum_{j=0}^{r} \binom{r}{j} (-1)^{r-j} \B{n+r-k}{j}.
\end{split}
\end{equation}

\vs

Let $k\geq r$. Then by (\ref{eqn_50}), we get
\begin{equation}\label{eqn_52}
\begin{split}
b_k^{(r)} & = \frac{1}{k!} \left\langle \left({e^{t}-1}\right)^r t^{k-r} \ | \  \B{n}{x}\right\rangle  = \frac{1}{k!} \left\langle \left({e^{t}-1}\right)^r  \ | \  t^{k-r} \B{n}{x}\right\rangle  \\
& = \frac{1}{k!}\binom{n}{k-r} (k-r)!\left\langle \left({e^{t}-1}\right)^r \ | \  \B{n+r-k}{x}\right\rangle  \\
 & = \frac{\binom{n}{k-r}}{r!\binom{k}{r}} \sum_{j=0}^{r} \binom{r}{j} (-1)^{r-j} \left\langle e^{jt} \ | \  \B{n+r-k}{x}\right\rangle  \\
 & = \frac{\binom{n}{k-r}}{r!\binom{k}{r}} \sum_{j=0}^{r} \binom{r}{j} (-1)^{r-j} \B{n+r-k}{j}.
\end{split}
\end{equation}
Therefore, by (\ref{eqn_49}), (\ref{eqn_51}) and (\ref{eqn_52}), we obtain the following theorem.
\begin{theorem}\label{thm_04}
For $n\in \iz_{+}$ and $r\in \nt$, we have
\end{theorem}
\begin{equation*}
\begin{split}
\B{n}{x} & = \sum_{k=0}^{r-1} \frac{\binom{r}{k}}{r!\binom{n+r-k}{r-k}} \sum_{j=0}^{r} \binom{r}{j} (-1)^{r-j} \B{n+r-k}{j} \Br{k}{r}{x} \\
& \quad + \sum_{k=r}^{n} \frac{\binom{n}{k-r}}{r!\binom{k}{r}} \sum_{j=0}^{r} \binom{r}{j} (-1)^{r-j} \B{n+r-k}{j} \Br{k}{r}{x}.
\end{split}
\end{equation*}

\section{Further Remarks}

For $n,m\in \iz_{+}$ with $n-m \geq 0$, we have
\begin{equation}\label{eqn_53}
\small
\begin{split}
\Br{n}{r}{x}&\Br{n-m}{r}{x} \\
& = \left(\sum_{l=0}^{n}\binom{n}{l} B_{n-l}^{(r-1)} \B{l}{x} \right) \left(\sum_{p=0}^{n-m}\binom{n-m}{p} B_{n-m-p}^{(r-1)} \B{p}{x} \right) \\
 & = \sum_{k=0}^{2n-m}\sum_{p=0}^{k}\binom{n-m}{p}\binom{n}{k-p} B_{n-m-p}^{(r-1)} B_{n-k+p}^{(r-1)} \B{p}{x}    \B{k-p}{x}.
\end{split}
\end{equation}

\vs

Let us consider $p(x) = E_n(x)$ with
\begin{equation}\label{eqn_54}
E_n(x) = p(x) = \sum_{k=0}^{n} b_k \B{k}{x}.
\end{equation}
Then we have
\begin{equation}\label{eqn_55}
p^{(k)}(x) = k!\binom{n}{k} E_{n-k}(x),
\end{equation}
and
\begin{equation}\label{eqn_56}
\begin{split}
b_k & = \frac{1}{k!} \left\langle \frac{e^{t}-1}{t} \ | \  p^{(k)}(x)\right\rangle  = \binom{n}{k} \left\langle \frac{e^{t}-1}{t} \ | \  E_{n-k}(x)\right\rangle  \\
& = \binom{n}{k} \frac{E_{n-k+1}(1) - E_{n-k+1}}{n-k+1} = -2 \binom{n}{k} \frac{E_{n-k+1}}{n-k+1}.
\end{split}
\end{equation}
By (\ref{eqn_54}) and (\ref{eqn_56}), we get
\begin{equation}\label{eqn_57}
E_n(x) = -2 \sum_{k=0}^{n} \binom{n}{k} \frac{E_{n-k+1}}{n-k+1} \B{k}{x}.
\end{equation}
From (\ref{eqn_57}), we can derive the following equation.
\begin{equation*}
\small
\begin{split}
E_n&(x)  E_{n-m}(x) \\
& = 4 \left(\sum_{l=0}^{n} \binom{n}{l} \frac{E_{n-l+1}}{n-l+1} \B{l}{x}\right) \left(\sum_{p=0}^{n-m} \binom{n-m}{p} \frac{E_{n-m-p+1}}{n-m-p+1} \B{p}{x}\right) \\
 & = 4 \sum_{k=0}^{2n-m} \sum_{l=0}^{k} \binom{n-m}{l} \binom{n}{k-l} \frac{E_{n-m-l+1}E_{n-k+l+1}}{(n-m-l+1)(n-k+l+1)} \B{l}{x} \B{k-l}{x}.
\end{split}
\end{equation*}
where $n,m\in \iz_{+}$ with $n-m\geq 0$.


\section*{Acknowledgements}

This research was supported by Basic Science Research Program through the National Research Foundation of Korea(NRF)
funded by the Ministry of Education, Science and Technology 2012R1A1A2003786.

\par\bigskip

\par

\bigskip\bigskip

\begin{center}\begin{large}

\end{large}\end{center}

\par

\end{document}